\documentclass[11pt,reqno,a4paper,oneside]{amsart}

\usepackage{latexsym}
\usepackage{amsmath}
\usepackage{amsthm}
\usepackage{amssymb}
\usepackage{amsfonts}
\usepackage{fancyhdr}
\usepackage{comment}

\usepackage{mathrsfs}                           

\newcommand{\mR}{\mathbf{R}}                    
\newcommand{\mC}{\mathbf{C}}                    
\newcommand{\mZ}{\mathbf{Z}}                    
\newcommand{\abs}[1]{\lvert #1 \rvert}          
\newcommand{\norm}[1]{\lVert #1 \rVert}         
\newcommand{\br}[1]{\langle #1 \rangle}         

\newcommand{\etilde}{\,\tilde{\rule{0pt}{6pt}}\,}

\newcommand{\supp}{\mathrm{supp}}

\newcommand{\curl}{\mathrm{curl}}

\newcommand{\closure}[1]{\overline{#1}}
\newcommand{\dbar}{\overline{\partial}}

\newcommand{\mOp}{\mathrm{Op}}

\theoremstyle{definition}
\newtheorem{thm}{Theorem}[section]
\newtheorem{prop}[thm]{Proposition}

\newtheorem{lemma}[thm]{Lemma}
\newtheorem*{definition}{Definition}

\newcommand{\ol}{\overline}
\newcommand{\Hscl}{H_{\text{scl}}}

\title[Determining Nonsmooth First Order Terms]{Determining Nonsmooth First Order Terms \\
from Partial Boundary Measurements}
\author{Kim
Knudsen$^1$ and Mikko Salo$^2$} \address{1\ Department of Mathematical
Sciences, Aalborg University}\address{2\ Department of Mathematics and Statistics / RNI,
University of Helsinki} \email{1\ kim\@@math.aau.dk} \email{2\
mikko.salo\@@helsinki.fi} \date{}

\date{}

\begin{document}

\begin{abstract}
  We extend results of Dos Santos Ferreira-Kenig-Sj\"ostrand-Uhlmann
  (arXiv:math.AP/0601466) to less smooth coefficients, and we show
  that measurements on part of the boundary for the magnetic
  Schr\"odinger operator determine uniquely the magnetic field related
  to a H\"older continuous potential. We
  give a similar result for determining a convection term. The
  proofs involve Carleman estimates, a smoothing procedure, and an
  extension of the Nakamura-Uhlmann pseudodifferential
  conjugation method to logarithmic Carleman weights.
\end{abstract}

\maketitle

\section{Introduction}

Let $\Omega \subseteq \mR^n$, $n \geq 3$, be a bounded simply connected $C^{\infty}$ domain with connected boundary. We will consider two inverse problems in $\Omega$. The first is the inverse problem for the magnetic Schr\"odinger operator, defined by 
\begin{equation*}
H_{W,V} = \sum_{j=1}^n (D_j + W_j)^2 + V,
\end{equation*}
where $D_j = \frac{1}{i} \frac{\partial}{\partial x_j}$, $W \in
L^{\infty}(\Omega ; \mC^n)$, and $V \in L^n(\Omega)$. Provided that
$0$ is not a Dirichlet eigenvalue of $H_{W,V}$ in $\Omega$, we define
the Dirichlet-to-Neumann map (DN map) formally as the magnetic normal
derivative
\begin{equation*}
N_{W,V}: f \mapsto (\nabla + iW)u \cdot \nu|_{\partial \Omega},
\end{equation*}
where $\nu$ is the outer unit normal to $\partial \Omega$, and $u \in
H^1(\Omega)$ solves the equation $H_{W,V} u = 0$ in $\Omega$ with
$u|_{\partial \Omega} = f$. Using a weak formulation, $N_{W,V}$ is
well defined as a bounded map from $H^{1/2}(\partial \Omega)$ into
$H^{-1/2}(\partial \Omega),$ see also Section~\ref{sec:prelim} below.

We consider $N_{W,V}$ as the boundary measurements related to the
operator $H_{W,V},$ and the inverse problem is to determine the coefficients of $H_{W,V}$ from partial
knowledge of $N_{W,V}.$ There is gauge equivalence: one has
$N_{W+\nabla p,V} = N_{W,V}$ when $p \in W^{1,\infty}(\Omega)$ and
$p|_{\partial \Omega} = 0$. Thus one may only expect to recover the
magnetic field $dW$ and the electric potential $V$ from the
boundary measurements. Here $dW$ is the $2$-form $d(\sum_{j=1}^n W_j(x)\,dx_j)$.

Let $x_0 \in \mR^n \smallsetminus \closure{\mathrm{ch}(\Omega)}$,
where $\mathrm{ch}(\Omega)$ is the convex hull of $\Omega$. We define
the front face of $\partial \Omega$ relative to $x_0$ by 
\begin{equation*}
F(x_0) = \{ x \in \partial \Omega \colon (x-x_0) \cdot \nu(x) \leq 0 \},
\end{equation*}
and we take $\tilde{F}$ to be an open neighborhood of $F(x_0)$ in
$\partial \Omega$. Also, let $W_{\mathrm{tan}} = W - (W \cdot \nu)\nu$
be the tangential component on $\partial \Omega$ of a vector field
$W$. 

Our main result states that the magnetic field and electric potential are determined by 
measurements $N_{W,V}$ made on the particular subset $\tilde{F}$.

\begin{thm} \label{thm:magnetic} Let $W_j \in C^{\varepsilon}(\ol
  \Omega ; \mC^n),\; \varepsilon > 0,$ and $V_j \in L^n(\Omega)$
  for $j = 1,2.$ Also assume that $0$ is not a Dirichlet eigenvalue of
  $H_{W_j,V_j}$ in $\Omega$. If
\begin{equation*}
N_{W_1,V_1} f|_{\tilde{F}} = N_{W_2,V_2} f|_{\tilde{F}} \quad \text{for all } f \in H^{1/2}(\partial \Omega),
\end{equation*}
then $dW_1 = dW_2$ in $\Omega$ and $(W_1)_{\mathrm{tan}} = (W_2)_{\mathrm{tan}}$ on $\partial \Omega$. If further $V_j \in L^{\infty}(\Omega)$, then $V_1 = V_2$ in $\Omega$.
\end{thm}

This theorem was proved for $W \in C^2(\ol \Omega)$ and $V \in
L^{\infty}(\Omega)$ in \cite{dksu}, following the $W = 0$ case in
\cite{ksu} (see also \cite{bukhgeimuhlmann}). These articles introduce
important ideas, and the main tools are Carleman estimates, the
construction of special solutions to the Schr\"odinger equation, and
analytic microlocal analysis. The special solutions generalize the
exponentially growing solutions in \cite{sylvesteruhlmann}.

We will employ the methods of \cite{dksu} and a smoothing argument
from \cite{paivarintapanchenkouhlmann} to work with $W
\in W^{1,n}(\Omega) \cap C^{\varepsilon}(\ol\Omega)$. The assumption
$W \in W^{1,n}(\Omega)$ ensures that the Carleman estimate
construction gives special solutions in $H^1(\Omega)$. To remove this
assumption, we will improve the construction of solutions by combining
Carleman estimates with the pseudodifferential conjugation method from
\cite{nakamurauhlmann}. The conjugation method was given in
\cite{nakamurauhlmann} for linear Carleman weights, and we extend this
method to logarithmic Carleman weights. This is a microlocal argument
using some ideas from \cite{duistermaathormander} and \cite{ksu}, and
it gives the theorem for $W \in C^{\varepsilon}(\ol\Omega)$. We also
allow complex coefficients and give a more precise proof when $\Omega$ has complicated geometry.

Earlier work on the magnetic inverse problem for $n \geq 3$ has dealt
with the case where the full DN map is known. The uniqueness question
was solved in \cite{sun} for small $W$ and in
\cite{nakamurasunuhlmann} for smooth $W$ without smallness
assumptions. The regularity of $W$ was relaxed to $C^1$ in
\cite{tolmasky} and to Dini continuous in \cite{salothesis}. A
constructive procedure for recovering $dW$ and $V$ from
$N_{W,V}$ is given in \cite{saloreconstruction}. The related inverse
scattering problem has been studied in \cite{eskinralston}. We also
mention \cite{brownsalo} and \cite{tzou} which consider boundary
determination and stability for magnetic Schr\"odinger operators, and
\cite{knudsen} which proves a partial data result for the nonsmooth
conductivity equation.

The other problem we consider is to determine a convection term from
boundary measurements. If $W \in L^{\infty}(\Omega ; \mR^n)$ is a
convection term and if $u \in H^1(\Omega)$ solves $(-\Delta + 2W \cdot
\nabla) u = 0$ in $\Omega$ with $u|_{\partial \Omega} = f$, we define
the related DN map formally as the normal derivative 
\begin{equation*}
N_W: f \mapsto \partial_{\nu} u|_{\partial \Omega}.
\end{equation*}
We remark that since $W$ is real valued, $0$ is not a Dirichlet
eigenvalue of $-\Delta + 2W \cdot \nabla$ by the maximum principle.
The map $N_W$ is again bounded from $H^{1/2}(\partial \Omega)$ into
$H^{-1/2}(\partial \Omega)$. The inverse problem of determining $W$
from $N_W$ can be reduced to the case of the magnetic Schr\"odinger
equation, and we obtain the following result as a corollary to Theorem
\ref{thm:magnetic}.

\begin{thm} \label{thm:convection}
  Let $W_j \in C^{\varepsilon}(\closure{\Omega} ; \mR^n)$,
  $\varepsilon > 0$, and $\nabla \cdot W_j \in L^{\infty}(\Omega)$ for
  $j = 1,2$. If
\begin{equation*}
N_{W_1} f|_{\tilde{F}} = N_{W_2} f|_{\tilde{F}} \quad \text{for all } f \in H^{1/2}(\partial \Omega),
\end{equation*}
then $W_1 = W_2$ in $\Omega$.
\end{thm}

Earlier work on the convection term problem for $n \geq 3$ includes
\cite{chengnakamurasomersalo} which solves the full data problem for
smooth $W$, and \cite{salothesis} which gives the same result for
Lipschitz continuous $W$.

The paper is organized as follows. In Section 2 we fix some notation
and preliminaries. Section 3 contains the facts on first order
elliptic equations which will be needed for the pseudodifferential
conjugation method. The Carleman estimates are in Section 4, and in
Section 5 we establish the pseudodifferential conjugation method for
logarithmic Carleman weights. Section 6 gives the construction of
special solutions to the magnetic Schr\"odinger equation, and Theorems
\ref{thm:magnetic} and \ref{thm:convection} are proved in Sections 7
and 8, respectively.

\section{Preliminaries}\label{sec:prelim}

First we fix some notation. We will always take $\Omega$ to be a
bounded simply connected domain in $\mR^n$, $n \geq 3$, with connected $C^{\infty}$
boundary. Recall that $\Omega$ is simply connected if it is connected and has trivial fundamental group (hence, the first de Rham cohomology vanishes). Let $\nu$ be the outer unit normal of $\Omega$, and
$\partial_{\nu} u = \nabla u \cdot \nu$ the normal derivative of $u$
on $\partial \Omega$.

We write $C(\closure{\Omega})$, $C^{\varepsilon}(\closure{\Omega})$,
and $C^{\infty}(\closure{\Omega})$ for the sets of complex valued
functions which are continuous, $\varepsilon$-H\"older continuous, and
infinitely differentiable, respectively, in $\closure{\Omega}$. The
corresponding spaces of compactly supported functions in $\Omega$ are
$C_c(\Omega)$, $C_c^{\varepsilon}(\Omega)$, and
$C^{\infty}_c(\Omega)$. The notation $C(\closure{\Omega} ; \mC^n)$
denotes the space of $n$-dimensional vector fields whose components
are in $C(\closure{\Omega})$, similarly for
$C^{\varepsilon}(\closure{\Omega} ; \mC^n),$ $L^p({\Omega} ; \mC^n)$
etc. For real-valued vector fields we write
$C^{\varepsilon}(\closure{\Omega} ; \mR^n)$ etc.

If $u, v \in L^2(\Omega)$ we write $(u|v) = \int_{\Omega} u \bar{v}
\,dx$, and also $(u|v) = \int_{\Omega} u \cdot \bar{v} \,dx$ when $u,
v \in L^2(\Omega ; \mC^n)$. The $L^2(\Omega)$ norm is denoted by $\norm{u} = (u|u)^{1/2}$. 
If $f, g \in L^2(\partial \Omega)$ let
$(f|g)_{\partial \Omega} = \int_{\partial \Omega} f \bar{g} \,dx$. We
will use the same notations for the pairing of distributions and
elements in the dual space. 

Denote by $H^s(\Omega)$ the $L^2$ Sobolev spaces in $\Omega$, and by
$W^{1,p}(\Omega)$ the usual $L^p$ Sobolev spaces. 
We will need the following well-known results concerning
multiplication in Sobolev spaces.

\begin{lemma} \label{lemma:multiplication}
If $a \in L^n(\Omega),\; n\geq 3,$ then $u \mapsto au$ maps $H^1(\Omega)$ to $L^2(\Omega)$ and 
\begin{align}
  \|au\|_{L^2(\Omega)} \leq C \|a\|_{L^n(\Omega)} \|u\|_{H^1(\Omega)}.\label{mult1}
\end{align}
If furthermore $a \in W^{1,n}\cap L^\infty(\Omega)$ then $u \mapsto au$
is bounded on $H^{s}(\Omega)$ for $-1 \leq s \leq 1$, and we have
\begin{align}
  \label{mult2}
  \|au\|_{H^{s}(\Omega)}\leq C \|a\|_{W^{1,n}\cap L^\infty(\Omega)}\|u\|_{H^{s}(\Omega)}.
\end{align}
\end{lemma}
\begin{proof}
  The first statement follows from the H\"older inequality and Sobolev
  embedding. The Leibniz rule and \eqref{mult1} give \eqref{mult2} for
  $s = 1$, the case $s = -1$ follows by duality, and the full result
  is obtained by interpolation.
\end{proof}

Next we note that when $W$ is H\"older continuous, one may always assume $\nabla \cdot W = 0$ after a gauge transformation which does not alter $N_{W,V}$ or the H\"older continuity of $W$.

\begin{lemma} \label{lemma:coulombgauge}
If $W \in C^{\varepsilon}(\closure{\Omega} ; \mC^n)$, then there exists $p \in C^{1+\varepsilon}(\closure{\Omega})$ such that $\nabla \cdot (W + \nabla p) = 0$ and $p|_{\partial \Omega} = 0$.
\end{lemma}
\begin{proof}
It is enough to choose $p$ to be the solution of $\Delta p = -\nabla \cdot W$ with $p|_{\partial \Omega} = 0$, and this Dirichlet problem has a $C^{1+\varepsilon}$ solution by \cite[Section 8.11]{gilbargtrudinger}.
\end{proof}

We proceed to give a precise definition of $N_{W,V}$. Let 
\begin{equation*}
H_{\Delta}(\Omega) = \{ u \in H^1(\Omega) \colon \Delta u \in L^2(\Omega) \}
\end{equation*}
with norm $\norm{u}_{H_{\Delta}(\Omega)} = \norm{u}_{H^1(\Omega)} +
\norm{\Delta u}_{L^2(\Omega)}$. This is a Banach space and
$C^{\infty}(\closure{\Omega})$ is a dense set. The trace operator is
bounded on $H_{\Delta}(\Omega)$ with values in
$H^{1/2}(\partial\Omega),$ and also $u \mapsto \partial_{\nu} u$ is bounded on
$H_{\Delta}(\Omega)$ with values in $H^{-1/2}(\partial\Omega)$. This follows by writing $u \in H_{\Delta}(\Omega)$ as $u = u_0 + u_1$ where $\Delta u_0 = 0$ with $u_0|_{\partial \Omega} = u|_{\partial \Omega}$, and $\Delta u_1 = \Delta u$ with $u_1|_{\partial \Omega} = 0$. If $\nabla \cdot W = 0$ then for $f \in H^{1/2}(\partial\Omega),$ the solution $u$ to
$H_{W,V} u = 0$ with boundary value $f$ is in $H_{\Delta}(\Omega)$. Therefore, the DN map $N_{W,V}: f \mapsto (\nabla + iW)u \cdot \nu|_{\partial \Omega}$ is bounded from $H^{1/2}(\partial \Omega)$ to $H^{-1/2}(\partial \Omega)$.


Finally, we will need the Green identity in the following form.
\begin{lemma} \label{lemma:hdelta_green_identity}
Let $W \in C(\closure{\Omega} ; \mC^n)$ with $\nabla \cdot W = 0$, and let $V \in L^n(\Omega)$. If $u \in H_{\Delta}(\Omega) \cap H^1_0(\Omega)$ and $v \in H_{\Delta}(\Omega)$, then 
\begin{equation*}
(H_{W,V}u|v) - (u|H_{\bar{W},\bar{V}}v) = -(\partial_{\nu} u|v)_{\partial \Omega}.
\end{equation*}
\end{lemma}
\begin{proof}
Let $u_j, v_j \in C^{\infty}(\closure{\Omega})$ with $u_j \to u$, $v_j \to v$ in $H_{\Delta}(\Omega)$. Then 
\begin{equation*}
(H_{W,V}u_j|v_j) - (u_j|H_{\bar{W},\bar{V}}v_j) = (u_j|\partial_{\nu} v_j)_{\partial \Omega} - (\partial_{\nu} u_j|v_j)_{\partial \Omega} -2i((W \cdot \nu)u_j|v_j)_{\partial \Omega}.
\end{equation*}
The claim follows by taking limits.
\end{proof}

\section{Elliptic equations of first order} \label{sec:firstorderelliptic}

To extend the Nakamura-Uhlmann pseudodifferential conjugation method to logarithmic Carleman weights, we will need to solve first order elliptic equations with variable coefficients in domains in $T^* \mR^n$. This is because the symbols of the conjugating operators arise as solutions to such equations. The main result in this section is Proposition \ref{prop:hmpsolution}, which shows the solvability of a first order equation related to the logarithmic Carleman weight.

The standard reference for the following facts is \cite{duistermaathormander}. See also \cite{ksu} for the specific case of limiting Carleman weights, which will be the case of interest for us. Thus, let $\varphi = \varphi(x)$ be a smooth real function in an open set $V \subseteq \mR^n$ with $\nabla \varphi \neq 0$ in $V$, and write $P = e^{\frac{\varphi}{h}}(-h^2\Delta)e^{-\frac{\varphi}{h}}$. The semiclassical Weyl symbol of $P$ is $p = a + ib$, where 
\begin{equation} \label{eq:abdef}
a(x,\xi) = \xi^2 - (\nabla \varphi)^2, \quad b(x,\xi) = 2\nabla \varphi \cdot \xi.
\end{equation}
We say that $\varphi$ is a limiting Carleman weight (for the Laplacian) if the Poisson bracket $\{a,b\} = \nabla_{\xi} a \cdot \nabla_x b - \nabla_x a \cdot \nabla_{\xi} b$ satisfies 
\begin{equation} \label{eq:limitingcarlemanpoisson}
\{a,b\} = 0 \quad \text{when } a = b = 0.
\end{equation}
This implies that $\{a,b\} = ca + db$ for some smooth $c,d$, and then on the set $\Sigma = \{a = b = 0\} \subseteq T^* V$ the Hamilton vector fields satisfy $[H_a,H_b] = cH_a + dH_b$. Recall that the Hamilton vector field of $f$ is $H_f = \nabla_{\xi} f \cdot \nabla_x - \nabla_x f \cdot \nabla_{\xi}$.

It follows that $\Sigma$ is an involutive manifold of codimension $2$, and the Frobenius theorem states that at each point of $\Sigma$ there are local coordinates in which $H_p = H_a + iH_b$ becomes an elliptic operator $a_1(y) \frac{\partial}{\partial y_1} + a_2(y) \frac{\partial}{\partial y_2}$. This implies local solvability of $H_p u = f$.

For our purposes local solvability is not enough, since we will need to solve a related equation in a full neighborhood of $\Sigma$. To do this we follow \cite{duistermaathormander} and find a smooth function $m$ satisfying $\{ma,mb\} = 0$ in a neighborhood of $\Sigma$, so that $H_{ma}$ and $H_{mb}$ commute near $\Sigma$ and in some new coordinates $H_{mp}$ becomes $\frac{\partial}{\partial y_1} + i \frac{\partial}{\partial y_2}$.

The first lemma states that in the present case where $a$ and $b$ come from a limiting Carleman weight, there is an explicit choice for $m$.

\begin{lemma}
If $m = \abs{\nabla \varphi(x)}^{-2}$, then $\{ma,mb\} = 0$ in $T^* V$.
\end{lemma}
\begin{proof}
From \eqref{eq:abdef} and \eqref{eq:limitingcarlemanpoisson} one gets (see \cite{ksu})
\begin{equation} \label{poissonbracketc1l1}
\{a,b\} = 4(\varphi''\xi \cdot \xi + \varphi''\nabla \varphi \cdot \nabla \varphi) = 4(c_1(x) a + (l_1(x) \cdot \xi)b)
\end{equation}
for some $c_1$ and $l_1$. Setting $\xi = 0$ and $\xi = \nabla \varphi$ gives $c_1 = -\frac{\varphi''\nabla \varphi \cdot \nabla \varphi}{(\nabla \varphi)^2}$ and $l_1 \cdot \nabla \varphi = -c_1$. By taking the terms in \eqref{poissonbracketc1l1} which are of second order in $\xi$ it follows that $(\varphi'' - c_1 I - 2l_1(\nabla \varphi)^t)\xi \cdot \xi = 0$, which implies 
\begin{equation*}
\varphi'' - c_1 I - l_1(\nabla\varphi)^t - \nabla\varphi(l_1)^t = 0.
\end{equation*}
Applying this matrix to $\nabla\varphi$ gives $l_1 =\frac{\varphi''\nabla \varphi}{(\nabla \varphi)^2}$.

Now $\nabla m = -2m l_1$, and 
\begin{align*}
\{ma,mb\} &= \nabla_{\xi}(ma) \cdot \nabla_x(mb) - \nabla_x(ma) \cdot \nabla_{\xi}(mb) \\
 &= m^2 \{a,b\} + mb \nabla_{\xi} a \cdot \nabla_x m - ma \nabla_x m \cdot \nabla_{\xi} b \\
 &= m^2 (\{a,b\} - 4bl_1 \cdot \xi + 4a l_1 \cdot \nabla\varphi) = 0.
\end{align*}
\end{proof}

We now specialize to the weight $\varphi(x) = \log\,\abs{x}$, and we will compute an explicit change of coordinates near $\Sigma$ which makes $H_{mp}$ into a $\dbar$ operator near $\Sigma$. In fact, the change of coordinates will consist of finding a codimension $2$ manifold in a neighborhood of $\Sigma$ which is transversal to the flows of $H_{ma}$ and $H_{mb}$. We then need to check that the flows originating from this manifold cover a full neighborhood of $\Sigma$ in the cotangent space.

Consider a truncated cone $V_0 = \{ x \in \mR^n \colon x_n > c \abs{x}, c < \abs{x} < c^{-1} \}$ for small $c > 0$. In suitable coordinates, $V_0$ will contain $\closure{\Omega}$. We will work in a fixed neighborhood 
\begin{equation*}
V = V(\delta) = \{ x \in \mR^n \colon x_n > (c-\delta) \abs{x}, c-\delta < \abs{x} < c^{-1}+\delta \}
\end{equation*}
with $\delta > 0$ small. Now $\Sigma$ is given by $\{(x,\xi) \in T^* V \colon x \cdot \xi = 0, \abs{x}\abs{\xi} = 1\}$, and we define a neighborhood 
\begin{equation*}
U = U(\delta) = \{(x,\xi) \in T^* V \colon \abs{x \cdot \xi} < \delta, \,1-\delta < \abs{x}\abs{\xi} < 1+\delta \}.
\end{equation*}
The first two new coordinates will be 
\begin{equation*}
y_1 = x \cdot \xi, \quad y_2 = \abs{x}\abs{\xi},
\end{equation*}
so $\Sigma$ is given by $\{ y_1 = 0, y_2 = 1 \}$.

We have $m = \abs{x}^2$ and $ma = \abs{x}^2 \abs{\xi}^2 - 1$, $mb = 2x \cdot \xi$, and the Hamilton vector fields are 
\begin{equation*}
H_{ma} = 2(\abs{x}^2 \xi \cdot \nabla_x - \abs{\xi}^2 x \cdot \nabla_{\xi}), \quad H_{mb} = 2(x \cdot \nabla_x - \xi \cdot \nabla_{\xi}).
\end{equation*}
These are smooth in $T^*\mR^n$, and it is possible to compute the flows explicitly. If $x$ and $\xi$ are linearly independent, then the flows starting from $(x,\xi)$ are given by 
\begin{align*}
\theta_{ma}(s,(x,\xi)) &= (e^{2y_1 s} \abs{x} ((\cos\,2\sqrt{y_2^2-y_1^2} s) \hat{x} + (\sin\,2\sqrt{y_2^2-y_1^2} s) J\hat{x}), \\
 &\hspace{19pt} e^{-2y_1 s} \abs{\xi} ((\sin\,2\sqrt{y_2^2-y_1^2} s) J\hat{\xi} + (\cos\,2\sqrt{y_2^2-y_1^2} s) \hat{\xi})), \\
\theta_{mb}(t,(x,\xi)) &= (e^{2t} x, e^{-2t} \xi),
\end{align*}
where $\hat{z} = \frac{z}{\abs{z}}$, and for $z$ in the oriented $2$-plane $T = \text{span}(x,\xi)$, $Jz$ is defined by 
\begin{equation*}
Jz = -(z \cdot \tilde{\xi})\hat{x} + (z \cdot \hat{x}) \tilde{\xi}
\end{equation*}
where $\tilde{\xi} = \frac{\xi-(\xi \cdot \hat{x})\hat{x}}{\abs{\xi-(\xi \cdot \hat{x})\hat{x}}}$. Thus, $Jz$ is the unique vector in $T$ for which $\abs{Jz} = \abs{z}$ and $(z,Jz)$ is a positive orthogonal basis of $T$.

The leaf $\Gamma_{x,\xi}$ through $(x,\xi) \in T^*V$ is the set of all points in $T^*V$ which can be reached from $(x,\xi)$ by the flows of $H_{ma}$ and $H_{mb}$.

\begin{lemma} \label{lemma:leafcharacterization}
If $(x,\xi) \in U$ then $\Gamma_{x,\xi} \subseteq U$ and 
\begin{equation*}
\Gamma_{x,\xi} = \{ (z,\eta) \in T^* V \colon z \cdot \eta = y_1, \ \abs{z} \abs{\eta} = y_2, \ \text{span}(z,\eta) = \text{span}(x,\xi) \}.
\end{equation*}
\end{lemma}
\begin{proof}
Denote the set on the right by $L_{x,\xi}$. It is easy to see that $\Gamma_{x,\xi}$ is contained in $L_{x,\xi}$ by checking that $z \cdot \eta$, $\abs{z} \abs{\eta}$, and $\text{span}(z,\eta)$ are constant along the flows. This also implies $\Gamma_{x,\xi} \subseteq U$.

For the converse, note that since $x \in V$ the plane $T = \text{span}(x,\xi)$ does not lie in $\{x_n = 0\}$, and we may define 
\begin{align*}
w &= w(x,\xi) = \frac{\text{proj}_T\,e_n}{\abs{\text{proj}_T\,e_n}}, \\
\zeta &= \zeta(x,\xi) = y_1 w + \sqrt{y_2^2-y_1^2} Jw.
\end{align*}
It follows that $(w, \zeta) \in L_{x,\xi}$. We claim that $(x,\xi) = \theta_{mb}(t,\theta_{ma}(s,(w,\zeta)))$, provided that 
\begin{align*}
s &= \frac{1}{2\sqrt{y_2^2 - y_1^2}} \angle (\hat{x} \cdot (w+iJw)), \\
t &= \frac{1}{2} \log\,\abs{x} - y_1 s,
\end{align*}
where we define $\angle(e^{i\alpha}) = \alpha$ for $-\pi < \alpha < \pi$. In fact we have 
\begin{equation*}
(\cos\,2\sqrt{y_2^2-y_1^2} s) w + (\sin\,2\sqrt{y_2^2-y_1^2} s) Jw = (\hat{x} \cdot w)w + (\hat{x} \cdot Jw)Jw = \hat{x}
\end{equation*}
and so $\theta_{ma}(s,(w,\zeta)) = (e^{2y_1 s} \hat{x}, \,\cdot\,)$ and $\theta_{mb}(t,\theta_{ma}(s,(w,\zeta))) = (x,\,\cdot\,)$. The claim follows since any two points in $L_{x,\xi}$ whose first components are the same must be identical. The construction also guarantees that the flow from $(w,\zeta)$ to $(x,\xi)$ stays in $T^* V$.

Now, if $(z,\eta) \in L_{x,\xi}$, then $w(z,\eta) = w(x,\xi)$ and $\zeta(z,\eta) = \zeta(x,\xi)$. This shows that one may reach both $(z,\eta)$ and $(x,\xi)$ by flows starting from the same point $(w,\zeta)$, which implies $(z,\eta) \in \Gamma_{x,\xi}$.
\end{proof}

Since the leaves are given by the points where $y_1$, $y_2$, and the plane $\text{span}(x,\xi)$ are constant, we take the next new variable to be this oriented plane. More precisely, we take $(y_3,\ldots,y_{2n-2})$ to be local coordinates corresponding to $\text{span}(x,\xi)$ on the Grassmannian $G(2,n)$, which consists of oriented two-planes in $\mR^n$. There is a single chart which achieves this, due to the fact that $\text{span}(x,\xi)$ does not lie in $\{x_n = 0\}$: one may apply an oriented version of Pl{\"u}cker coordinates, or if $n = 3$ it is sufficient to use the identification $G(2,3) = S^2$ and stereographic projection.

The final coordinates will be the flow variables $y_{2n-1} = s$ and $y_{2n} = t$, where $s$ and $t$ were given in the proof of Lemma \ref{lemma:leafcharacterization}. The codimension $2$ manifold transversal to the flows will then be $\{ y_{2n-1} = y_{2n} = 0 \}$. We have arrived at the desired change of coordinates.

\begin{lemma}
The map $\Phi: (x,\xi) \mapsto y$ defined in the discussion above is smooth and injective on $U(\delta)$ for $\delta$ small, and it is a diffeomorphism onto its image in $\mR^{2n}$. In the new coordinates $H_{mp}$ becomes $\frac{\partial}{\partial s} + i \frac{\partial}{\partial t}$.
\end{lemma}
\begin{proof}
If $\Phi(x,\xi) = \Phi(x',\xi')$ where $(x,\xi), (x',\xi') \in U$, then Lemma \ref{lemma:leafcharacterization} implies that $(x',\xi')$ is on the leaf through $(x,\xi)$. Also, the vectors $w, \zeta$ and the coordinates $s, t$ in Lemma \ref{lemma:leafcharacterization} are the same whether they are computed from $(x,\xi)$ or $(x',\xi')$. It follows that 
$(x,\xi) = \theta_{mb}(t,\theta_{ma}(s,(w,\zeta))) = (x',\xi')$, which shows that $\Phi$ is injective. Since $\Phi$ is smooth, it is a diffeomorphism onto its image. Further, $H_{mp}$ becomes $\frac{\partial}{\partial s} + i \frac{\partial}{\partial t}$ because the flows of $H_{ma}$ and $H_{mb}$ commute.
\end{proof}

Note that $H_{mp} = mH_p$ on $\Sigma$, so the following result states in particular that one may solve $H_p u = f$ on $\Sigma$.

\begin{prop} \label{prop:hmpsolution}
If $f(x,\xi) \in C^{\infty}(\closure{U(2\delta)})$ for $\delta$ small, then the equation 
\begin{equation*}
H_{mp} u = mf \quad \text{in } U(\delta)
\end{equation*}
has a solution $u \in C^{\infty}(U(\delta))$ satisfying for all $k > 0$ 
\begin{equation*}
\abs{\partial_x^{\alpha} \partial_{\xi}^{\beta} u(x,\xi)} \leq C_k \norm{f}_{W^{k,\infty}(U(2\delta))}, \quad (x,\xi) \in U(\delta)
\end{equation*}
whenever $\abs{\alpha} + \abs{\beta} \leq k$.
\end{prop}
\begin{proof}
Writing $u = \tilde{u} \circ \Phi$, it is enough to solve 
\begin{equation*}
(\partial_{y_{2n-1}} + i\partial_{y_{2n}})\tilde{u} = (\chi m f) \circ \Phi^{-1} \quad \text{on } \Phi(U(2\delta))
\end{equation*}
where $\chi \in C_c^{\infty}(U(2\delta))$ and $\chi = 1$ near $U(\delta)$. This may be solved using the Cauchy transform $(\partial_{y_{2n-1}} + i\partial_{y_{2n}})^{-1}$, and the norm estimates are an immediate consequence.
\end{proof}

\section{Carleman estimates}

In this section we will first recall the Carleman estimates for
$H_{W,V}$ in \cite{dksu}, and we note that the shifted estimate is
valid when $W \in W^{1,n} \cap L^{\infty}(\Omega ; \mC^n)$. We then
use the Carleman estimate to solve an equation involving a conjugated
version of $H_{W,V}$.

The Carleman estimate in \cite{dksu} is an estimate for the magnetic
Schr\"odinger equation proved for limiting Carleman weights. Denote by
$\tilde \Omega\subseteq \mR^n$ an open set such that $\ol \Omega \subseteq
\tilde \Omega$, and recall that $\varphi$ is a limiting Carleman
weight in $\tilde{\Omega}$ if \eqref{eq:abdef} and
\eqref{eq:limitingcarlemanpoisson} are satisfied in $T^*
\tilde{\Omega}$. Later we will restrict ourselves to the particular
logarithmic weight. Below we write $\partial \Omega_{\pm} = \{ x \in \partial \Omega \colon \pm \partial_{\nu} \varphi(x) \geq 0\}$, and $A \lesssim B$ if $A \leq CB$ where $C > 0$ is a constant independent of $h$ and $\varepsilon$. Recall that $\norm{\,\cdot\,} = \norm{\,\cdot\,}_{L^2(\Omega)}$.

\begin{prop}\label{prop:carleman}
  Let $\varphi$ be a limiting Carleman weight on $\tilde \Omega.$ Suppose
  $W \in W^{1,n}\cap L^\infty(\Omega;\mC^n)$ and $V\in L^n(\Omega).$
  Then for $u\in H_{\Delta} \cap H^1_0(\Omega)$ and $h$ small, we have the Carleman
  estimate
  \begin{multline}\label{eq:carleman}
    -h(\partial_{\nu} \varphi e^{\frac{\varphi}{h}} \partial_{\nu} u| e^{\frac{\varphi}{h}} \partial_{\nu} u)_{\partial \Omega_{-}} + \|e^{\frac{\varphi}{h}}u\|^2 + \|e^{\frac{\varphi}{h}}h\nabla u\|^2 
     \\ \lesssim
    h^2\|e^{\frac{\varphi}{h}} H_{W,V}u\|^2 +
    h(\partial_{\nu} \varphi e^{\frac{\varphi}{h}} \partial_{\nu} u| e^{\frac{\varphi}{h}} \partial_{\nu} u)_{\partial \Omega_{+}}.
  \end{multline}
\end{prop}
\begin{proof}
  By introducing $v =e^{\frac{\varphi}{h}}u$ it follows that \eqref{eq:carleman}
  is equivalent to the a priori estimate
  \begin{multline}\label{eq:carleman2}
    -h(\partial_{\nu} \varphi \partial_{\nu} v | \partial_{\nu} v)_{\partial \Omega_{-}} + \|v\|^2 + \|h\nabla v\|^2 \\
     \lesssim
    h^2\|e^{\frac{\varphi}{h}} H_{W,V}e^{-\frac{\varphi}{h}} v\|^2 +
    h(\partial_{\nu} \varphi \partial_{\nu} v | \partial_{\nu} v)_{\partial \Omega_{+}}
  \end{multline}
  for the conjugated operator $e^{\frac{\varphi}{h}} H_{W,V}e^{-\frac{\varphi}{h}}.$ 
  To prove \eqref{eq:carleman2} the idea is to work with the convexified weight
  \begin{align*}
    \tilde \varphi (x) = \varphi(x) + h\frac{\varphi(x)^2}{2\varepsilon}
  \end{align*}
  and first obtain the estimate for $e^{\frac{\tilde{\varphi}}{h}}
  H_{W,V}e^{-\frac{\tilde{\varphi}}{h}}.$ This operator can be split into a self-adjoint term $\tilde P$, a
  skew-adjoint term $i\tilde Q$, and a remainder term $\tilde R$, 
  \begin{align*}
    e^{\frac{\tilde{\varphi}}{h}} H_{W,V}e^{-\frac{\tilde{\varphi}}{h}}  &= h^{-2}(\tilde P +
    i\tilde Q + \tilde R),
  \end{align*}
  where 
  \begin{align*}
    \tilde P &= -h^2\Delta - (\nabla \tilde\varphi)^2, \\
    \tilde Q &= hD \circ \nabla \tilde\varphi
    +\nabla \tilde\varphi \circ hD, \\
    \tilde R &= 2hW \cdot (hD + i\nabla\tilde\varphi)+ h^2(W^2 + D \cdot W +V).
  \end{align*}
  We then have from \cite[Eq.~(2.12)]{dksu} the estimate
  \begin{align}\label{eq:carlemanfree}
    \frac{h^2}\varepsilon(\|v\|^2 + \|h\nabla v\|^2) \lesssim
   \|( \tilde P + i \tilde Q)v\|^2 +
    h^3(\partial_\nu\tilde\varphi \partial_\nu v | \partial_\nu v)_{\partial\Omega}
  \end{align}
  for some sufficiently small $\varepsilon>0.$  By the assumptions on $W$ and $V$, it
  follows from \eqref{mult1} that 
  \begin{align*}
    \|\tilde Rv\|^2 &\lesssim h^2  (\|W\|_{L^\infty(\Omega)}^2  + 
    \|W^2+D\cdot W + V\|^2 _{L^n(\Omega)})(\|v\|^2  + \|h\nabla v\|^2 ) \\
    &\lesssim h^2 (\|v\|^2 +\|h\nabla v\|^2 ).
  \end{align*}
  Hence by choosing $\varepsilon>0$ sufficiently small,
  the term $\|\tilde Rv\|^2$ can be absorbed in the left hand side of
  \eqref{eq:carlemanfree} and consequently
  \begin{align*}
    \frac{h^2}\varepsilon(\|v\|^2 + \|h\nabla v\|^2) \lesssim
   h^4\| e^{\frac{\tilde{\varphi}}{h}} H_{W,V}e^{-\frac{\tilde{\varphi}}{h}} v \|^2 +
    h^3(\partial_\nu\tilde\varphi \partial_\nu v | \partial_\nu v)_{\partial\Omega},
  \end{align*}
  from which \eqref{eq:carleman2} can be derived as in \cite{dksu}.
\end{proof}
Next we show that the Carleman estimate \eqref{eq:carleman} can be shifted to a lower Sobolev index. For this we will use the semiclassical Sobolev spaces $H^{s}_{\text{scl}}(\mR^n)$ with norm $\norm{f}_{H^s_{\text{scl}}(\mR^n)} = \norm{\br{hD}^s f}_{L^2(\mR^n)}$, where $\br{\xi} = (1+\abs{\xi}^2)^{1/2}$.
\begin{prop}
  Let $\varphi$ be a limiting Carleman weight, and suppose that $W\in W^{1,n}\cap
  L^\infty(\Omega;\mC^n)$ and $V\in L^n(\Omega).$ Then for $v\in C^\infty_c(\Omega)$  we have the estimate
  \begin{align}\label{carleman2}
    \|v\|_{L^2(\mR^n)}^2 \lesssim h^2\|e^{\frac{\varphi}{h}} H_{W,V}
    e^{-\frac{\varphi}{h}}u\|_{\Hscl^{-1}(\mR^n)}^2.
  \end{align}
\end{prop}
\begin{proof}
  Using the notation in Proposition \ref{prop:carleman} we have from
  \cite[Proposition 2.4]{dksu} the estimate
  \begin{align*}
    \frac{h^2}\varepsilon\|v\|^2_{H^1_{\text{scl}}(\mR^n)} \lesssim \|(\tilde P + i\tilde Q)\br{hD}v\|^2_{\Hscl^{-1}(\mR^n)},
  \end{align*}
  for $v\in C_c^\infty(\hat \Omega),$ where $\hat\Omega$ is open and
  $\ol\Omega\subseteq\hat\Omega$ and $\closure{\hat{\Omega}} \subseteq \tilde \Omega.$ Suppose we have extended
  $W,V$ to $\mR^n.$ The semiclassical counterparts of \eqref{mult1} and \eqref{mult2} are 
  \begin{align*}
   h\norm{au}_{L^2(\mR^n)} &\lesssim \norm{a}_{L^n(\mR^n)} \norm{u}_{H^1_{\text{scl}}(\mR^n)}, \\
   \norm{au}_{H^{1}_{\text{scl}}(\mR^n)} &\lesssim \norm{a}_{W^{1,n} \cap L^{\infty}(\mR^n)} \norm{u}_{H^{1}_{\mathrm{scl}}(\mR^n)},
  \end{align*}
  and the corresponding dual estimates yield 
  \begin{equation*}
    \norm{\tilde{R}w}_{H^{-1}_{\text{scl}}(\mR^n)} \lesssim h \norm{w}_{L^2(\mR^n)}.
  \end{equation*}
  Hence for $\varepsilon>0$ sufficiently small we have
  \begin{align*}
    \frac{h^2}\varepsilon\|v\|^2_{H^1_{\text{scl}}(\mR^n)} \lesssim
    h^4 \| e^{\frac{\tilde{\varphi}}{h}} H_{W,V}
    e^{-\frac{\tilde{\varphi}}{h}} \br{hD} v\|^2_{\Hscl^{-1}(\mR^n)},
  \end{align*}
  from which \eqref{carleman2} can be derived as in \cite{dksu}.
\end{proof}

As a consequence of the estimate \eqref{carleman2} and the Hahn-Banach
theorem, we have the following solvability result.

\begin{prop} \label{prop:firstinhomog}
  Let $\varphi$ be a limiting Carleman weight, and let $W\in W^{1,n}\cap
  L^\infty(\Omega;\mC^n)$ and $V\in L^n(\Omega).$ If $h$ is sufficiently small, then for any $f\in L^2(\Omega)$ the equation
  \begin{align} \label{inhomogeq}
    e^{\frac{\varphi}{h}} H_{W,V} e^{-\frac{\varphi}{h}}u = f
  \end{align}
has a solution $u \in H^1(\Omega)$ with 
\begin{align} \label{inhomognormestimates}
  \|u\|_{L^2(\Omega)} + \|hDu\|_{L^2(\Omega)} \lesssim h\|f\|_{L^2(\Omega)}.
 \end{align}
\end{prop}

\section{Pseudodifferential conjugation}

In this section we will prove Proposition \ref{prop:firstinhomog} in
the case of uniformly continuous coefficients. This will follow by
extending the pseudodifferential conjugation technique introduced in
\cite{nakamurauhlmann} (see also \cite{nakamurauhlmannerratum}) to
logarithmic Carleman weights
\begin{align}
  \label{phidef}
  \varphi(x) = \log\,\abs{x-x_0},\quad x_0 \notin \closure{\Omega}.
\end{align}

The main result is the following.

\begin{prop} \label{prop:secondinhomog}
Let $\varphi(x)$ be as in \eqref{phidef}, and assume that $W \in C(\closure{\Omega} ; \mC^n)$ and $\nabla \cdot W,\, V \in L^n(\Omega ; \mC)$. If $h$ is sufficiently small, then for any $f \in L^2(\Omega)$ the equation \eqref{inhomogeq} 
has a solution $u \in H^1(\Omega)$ satisfying \eqref{inhomognormestimates}.
\end{prop}

To prove this, we start by choosing $\tilde{\Omega}$ to be an open set containing $\closure{\Omega}$ with $x_0 \notin \closure{\tilde{\Omega}}$, and we take coordinates so that $x_0 = 0$ and $\closure{\tilde{\Omega}} \subseteq \{ x_n > c\abs{x}, c < \abs{x} < c^{-1} \}$ for $c$ small. Note that $H_{W,V} = -\Delta + 2W \cdot D + \tilde{V}$ with $\tilde{V} \in L^n(\tilde\Omega)$. We ignore the zero order term for the moment and consider the conjugated operator 
\begin{equation*}
e^{\frac{\varphi}{h}} ((hD)^2 + 2h W \cdot hD) e^{-\frac{\varphi}{h}} = P + hQ
\end{equation*}
where $P = (hD+i\nabla\varphi)^2$ and $Q = 2W \cdot (hD + i\nabla\varphi)$.

If $h$ is sufficiently small, Proposition \ref{prop:firstinhomog} (in the case $W = V = 0$) implies that we may solve $Pu = f$ in $\tilde{\Omega}$ for any $f \in L^2(\tilde{\Omega})$, and the solution operator $P^{-1}: f \mapsto u$ is a linear map which satisfies 
\begin{equation} \label{pinvestimates}
\norm{P^{-1}f}_{L^2(\tilde{\Omega})} + \norm{hD P^{-1}f}_{L^2(\tilde{\Omega})} \lesssim h^{-1} \norm{f}_{L^2(\tilde{\Omega})}.
\end{equation}
We would like to solve $(P+hQ)u = f$ in the same way. If
$\norm{W}_{L^{\infty}(\tilde \Omega)}$ is small then $P+hQ$ is a small perturbation of $P$, and a solution is obtained from the Neumann series. If $W$ is large a different method is needed. In this case we will find order $0$ pseudodifferential operators which conjugate the operator into a small perturbation of $P$. We use standard classes of semiclassical symbols and operators.

\begin{definition}
If $0 \leq \sigma < 1/2$ and $m \in \mR$, let $S^m_{\sigma}$ be the space of all functions $a(x,\xi;h)$ where $x,\xi \in \mR^n$ and $h \in (0,h_0]$, $h_0 \leq 1$, such that $a$ is smooth in $x$ and $\xi$ and 
\begin{equation*}
\abs{\partial_x^{\alpha} \partial_{\xi}^{\beta} a(x,\xi;h)} \leq C_{\alpha \beta} h^{-\sigma\abs{\alpha+\beta}} (1+\abs{\xi}^2)^{m/2}
\end{equation*}
for all $\alpha,\beta$. If $a \in S^m_{\sigma}$ we define an operator $A = \mOp_h(a)$ by 
\begin{equation*}
Af(x) = (2\pi)^{-n} \iint_{\mR^{2n}} e^{i(x-y)\cdot\xi} a\left(\frac{x+y}{2},h\xi\right) f(y) \,d\xi \,dy.
\end{equation*}
\end{definition}

Note that we define the operators using semiclassical Weyl quantization. Operators in $S^0_{\sigma}$ are bounded on $L^2$ with norm uniformly bounded in $h$, and the composition of two operators is again an operator in the same class. We assume familiarity with semiclassical calculus in what follows, for more details see \cite{dimassisjostrand} and also \cite{saloreconstruction} where the result of this section was proved in the semiclassical setup for linear Carleman weights.

To deal with the nonsmooth coefficients we extend $W$ to a vector field in $C_c(\tilde{\Omega} ; \mC^n)$ and $\tilde{V}$ by zero, and we consider the mollifications 
\begin{equation} \label{wapproximation}
W^{\sharp} = W \ast \chi_{\delta}, \quad \tilde{V}^{\sharp} = \tilde{V} \ast \chi_{\delta}
\end{equation}
where $\delta = h^{\sigma}$ and $0 < \sigma < 1/2$. Here $\chi_{\delta}(x) = \delta^{-n} \chi(x/\delta)$ is the usual mollifier with $\chi \in C_c^{\infty}(\mR^n)$, $0 \leq \chi \leq 1$, and $\int \chi \,dx = 1$. We write $W^{\flat} = W - W^{\sharp}$ and $\tilde{V}^{\flat} = \tilde{V} - \tilde{V}^{\sharp}$ and note the estimates 
\begin{align}
\norm{\partial^{\alpha} W^{\sharp}}_{L^{\infty}}, \norm{\partial^{\alpha} \tilde{V}^{\sharp}}_{L^n} &= O(h^{-\sigma \abs{\alpha}}), \label{wvsharpest} \\
\norm{W^{\flat}}_{L^{\infty}}, \norm{\tilde{V}^{\flat}}_{L^n} &= o(1) \label{wvflatest}
\end{align}
as $h \to 0$.

Write $q^{\sharp} = 2W^{\sharp} \cdot (\xi + i\nabla\varphi) - hD
\cdot W^{\sharp}$, so that $Q^{\sharp} = \mOp_h(q^{\sharp}) = 2 W^{\sharp} \cdot (hD + i\nabla\varphi)$. The following lemma shows that one may use pseudodifferential operators to conjugate away the main part $Q^{\sharp}$ of the first order term $Q$.

\begin{lemma} \label{lemma:psdoconjugation}
There exist $c, \tilde{c}, r \in S^0_{\sigma}$ so that 
\begin{equation*}
(P+hQ^{\sharp})C = \tilde{C}P + h^{2-2\sigma} R \quad \text{in } \Omega.
\end{equation*}
Both $C$ and $\tilde{C}$ are elliptic, in the sense that $c$ and $\tilde{c}$ are nonvanishing for small $h$.
\end{lemma}
\begin{proof}
Suppose $c$ is any symbol in $S^0_{\sigma}$ which is equal to $1$ outside a large ball in $(x,\xi)$. Then for some $r_0 \in S^0_{\sigma}$, 
\begin{equation} \label{conjugationcommute}
(P+hQ^{\sharp})C = CP + h \mOp_h(\frac{1}{i} H_p c + q^{\sharp}c) + h^{2-2\sigma} R_0 \quad \text{in } \Omega.
\end{equation}
We will choose $c$ so that $\frac{1}{i} H_p c + q^{\sharp}c$ is of lower order. Trying $c = e^{i\phi}$ and using the notation in Section \ref{sec:firstorderelliptic} and Proposition \ref{prop:hmpsolution}, we choose $\phi$ to be a solution of 
\begin{equation} \label{phiequation}
H_{mp}\phi = -mq^{\sharp} \quad \text{in } U = U(\delta).
\end{equation}
Since $x$ and $\xi$ are bounded in $U$, \eqref{wvsharpest} gives 
\begin{equation*}
\abs{\partial_x^{\alpha} \partial_{\xi}^{\beta} \phi(x,\xi)} \leq C_{\alpha \beta} h^{-\sigma\abs{\alpha+\beta}}, \quad (x,\xi) \in U.
\end{equation*}
Let $\chi(x,\xi)$ be a smooth cutoff supported in $U(\delta)$ which is equal to $1$ on $U(\delta/2)$. Then $c$ is chosen to be 
\begin{equation*}
c = e^{i\chi\phi}.
\end{equation*}
It follows that $c \in S^0_{\sigma}$ and $c$ is nonvanishing.

Since $H_{mp} = mH_p + pH_m$, we get by \eqref{phiequation} 
\begin{align*}
m(\frac{1}{i} H_p c + q^{\sharp} c) &= ((H_{mp} \chi)\phi + \chi H_{mp} \phi - p H_m(\chi\phi) + mq^{\sharp})c \\
 &= ((H_{mp} \chi)\phi + (1-\chi)mq^{\sharp} -p H_m(\chi\phi))c.
\end{align*}
Let $\psi \in C^{\infty}_c(\tilde{\Omega})$ and $\psi = 1$ near $\closure{\Omega}$. Then $\psi(\frac{1}{i} H_p c + q^{\sharp} c) = lp$ where 
\begin{equation*}
l = \psi \Big( \frac{1-\chi}{p} q^{\sharp} + \frac{1}{mp} (H_{mp} \chi)\phi - \frac{1}{m} H_m(\chi\phi) \Big)c
\end{equation*}
and $l \in h^{-\sigma} S^0_{\sigma}$. Then \eqref{conjugationcommute} gives 
\begin{equation*}
(P+hQ^{\sharp})C = \tilde{C}P + h^{2-2\sigma} R \quad \text{in } \Omega
\end{equation*}
where $\tilde{c} = c + hl$ is of order $0$ and elliptic, and $r \in S^0_{\sigma}$.
\end{proof}

We may now prove the main result. The proof involves a number of cutoffs to restrict the functions to the set $\tilde{\Omega}$ where $\varphi$ is smooth.

\begin{proof}
(of Proposition \ref{prop:secondinhomog}) Using the notation in this section, we want to solve 
\begin{equation} \label{rewritteneq}
(P + hQ^{\sharp} + hQ^{\flat} + h^2\tilde{V}^{\sharp} + h^2\tilde{V}^{\flat})u = h^2 f \quad \text{in } \Omega
\end{equation}
where $Q^{\flat} = 2 W^{\flat} \cdot (hD + i\nabla\varphi)$. Let $\psi_j$ ($j \geq 1$) be $C^{\infty}_c(\tilde{\Omega})$ functions with $\psi_1 = 1$ near $\closure{\Omega}$ and $\psi_{j+1} = 1$ near $\supp(\psi_j)$. We try a solution 
\begin{equation} \label{uansatz}
u = C \psi_2 P^{-1} \psi_3 \tilde{C}^{-1} \psi_4 v
\end{equation}
for some $v \in L^2(\mR^n)$. Here $\tilde{C}^{-1}$ is the inverse of $\tilde{C}$ on $L^2(\mR^n)$, which exists for small $h$. Inserting \eqref{uansatz} in \eqref{rewritteneq} and using Lemma \ref{lemma:psdoconjugation} gives 
\begin{equation} \label{mplusr1}
(M+R_1)v = h^2 f \quad \text{in } \Omega
\end{equation}
where $M = \tilde{C} P \psi_2 P^{-1} \psi_3 \tilde{C}^{-1} \psi_4$ and $\norm{\psi_1 R_1}_{L^2 \to L^2} = o(1)$ as $h \to 0$. Here we have used \eqref{mult1}, \eqref{pinvestimates}, \eqref{wvsharpest}, \eqref{wvflatest}, and the identity $\partial_j C = C \partial_j + \mOp_h(\partial_{x_j} c)$.

We may write 
\begin{equation*}
M = \psi_2 + [\tilde{C},\psi_2] \tilde{C}^{-1} \psi_4 + \tilde{C} [P,\psi_2] P^{-1} \psi_3 \tilde{C}^{-1} \psi_4.
\end{equation*}
Extending $f$ by zero and modifying \eqref{mplusr1} slightly, we obtain an equation in $\mR^n$ of the form 
\begin{equation} \label{iplusm1plusr1}
(I+M_1+\psi_1 R_1)v = h^2 f \quad \text{in } \mR^n
\end{equation}
where $M_1 = \psi_1 [\tilde{C},\psi_2] \tilde{C}^{-1} \psi_4 + \psi_1 \tilde{C} [P,\psi_2] P^{-1} \psi_3 \tilde{C}^{-1} \psi_4$. By the pseudolocal property $\norm{M_1}_{L^2 \to L^2} = o(1)$ as $h \to 0$, so we obtain a solution $v \in L^2(\mR^n)$ to \eqref{iplusm1plusr1} by Neumann series. Since $\psi_j = 1$ on $\closure{\Omega}$ the function $v$ will satisfy \eqref{mplusr1} near $\closure{\Omega}$. It follows that \eqref{uansatz} solves \eqref{rewritteneq}, and one obtains the desired norm estimates from \eqref{pinvestimates}.
\end{proof}

\section{Construction of solutions} \label{sec:construction}

We now give the construction of special solutions to $H_{W,V} u = 0$ in $\Omega$. This proceeds as in \cite{dksu}, except that an additional smoothing argument is required. The solutions will have the form $u = e^{-\frac{\varphi}{h}} v$ where $v$ is a WKB solution for the conjugated operator $e^{\frac{\varphi}{h}} H_{W,V} e^{-\frac{\varphi}{h}}$. Thus, $u$ is of the form 
\begin{equation} \label{WKBansatz}
u = e^{-\frac{\rho}{h}}(a+r)
\end{equation}
where $\rho = \varphi + i\psi$ and $\psi$ is a real valued phase function, $a$ is an amplitude, and $r$ is a correction term. By Lemma \ref{lemma:coulombgauge} we may assume that $\nabla \cdot W = 0$, and then 
\begin{equation*}
h^2 H_{W,V} = (hD)^2 + 2hW \cdot hD + h^2 \tilde{V}
\end{equation*}
with $\tilde{V} \in L^n$, and inserting \eqref{WKBansatz} into $h^2 H_{W,V} u = 0$ gives the equation 
\begin{equation*}
((hD+i\nabla\rho)^2 + 2hW \cdot (hD+i\nabla\rho) + h^2 \tilde{V})(a+r) = 0.
\end{equation*}
Collecting like powers of $h$, $u$ will be a solution provided that one has in $\Omega$ 
\begin{align}
& (\nabla\rho)^2 = 0, \label{eikonaleq} \\
& (\nabla\rho \cdot D + \nabla\rho \cdot W + \frac{1}{2i} \Delta\rho)a = 0, \label{transporteq} \\
& e^{\frac{\rho}{h}} h^2 H_{W,V} e^{-\frac{\rho}{h}} r = -h^2 H_{W,V} a. \label{correctioneq}
\end{align}

We fix $\varphi(x)$ defined by \eqref{phidef}. Now \eqref{eikonaleq} is an eikonal equation for $\psi$, which reads 
\begin{equation*}
(\nabla \psi)^2 = (\nabla \varphi)^2, \qquad \nabla \varphi \cdot \nabla \psi = 0.
\end{equation*}
As shown in \cite{ksu}, a solution is given by 
\begin{equation} \label{psidef}
\psi(x) = \mathrm{dist}_{S^{n-1}} \Big( \frac{x-x_0}{\abs{x-x_0}},\omega \Big)
\end{equation}
where $\omega \in S^{n-1}$ is chosen so that $\psi$ is smooth in $\tilde{\Omega}$.

We choose coordinates so that $x_0 = 0$ and $\closure{\tilde{\Omega}} \subseteq \{ x_n > 0 \}$, and we take $\omega = e_1$.
We also write $x = (x_1,r\theta)$ with $r > 0$ and $\theta \in S^{n-2}$, and consider the change of variables $\Psi: x \mapsto (z,\theta)$ where $z = x_1 + ir$ is a complex variable. Writing $\tilde{f} = f \circ \Psi^{-1}$, we get $\tilde{\rho} = \log\,z$, $(\nabla \rho)\etilde = \frac{1}{z}(e_1+ie_r)$ where $e_r = (0,\theta)$, $(\Delta \rho)\etilde = -\frac{2(n-2)}{z(z-\bar{z})}$, and $\nabla \rho \cdot \nabla$ becomes $\frac{2}{z} \partial_{\bar{z}}$.

We see that \eqref{transporteq} is a $\dbar$-equation in the new coordinates, and the solution $a$ will have in general the same regularity as $W$. On the other hand, in \eqref{correctioneq} one needs two derivatives of $a$ on the right. To deal with this for $W \in C^{\varepsilon}(\closure{\Omega}; \mC^n)$ we use the same approximation procedure as in \eqref{wapproximation}, now with $\sigma > 0$ small and $W$ extended as a vector field in $C_c^{\varepsilon}(\tilde{\Omega} ; \mC^n)$. The solution will be taken of the form $u = e^{-\frac{\rho}{h}}(a+r)$, where $\rho$ is as above and one has in $\Omega$ 
\begin{align}
& (\nabla\rho \cdot D + \nabla\rho \cdot W^{\sharp} + \frac{1}{2i} \Delta\rho)a = 0, \label{transporteq_mollified} \\
& e^{\frac{\rho}{h}} h^2 H_{W,V} e^{-\frac{\rho}{h}} r = - 2ih (\nabla\rho \cdot W^{\flat})a - h^2 H_{W,V} a. \label{correctioneq_mollified}
\end{align}
To solve \eqref{transporteq_mollified} we take $a = \tilde{a} \circ \Psi$, so $\tilde{a}$ must satisfy 
\begin{equation*}
\Big(\partial_{\bar{z}} + \frac{i}{2}(e_1+i e_r) \cdot \tilde{W}^{\sharp} - \frac{n-2}{2(z-\bar{z})}\Big) \tilde{a} = 0 \quad \text{in } \Psi(\Omega).
\end{equation*}
Here $\tilde{W}^{\sharp} = W^{\sharp} \circ \Psi^{-1}$. We try $\tilde{a} = (z-\bar{z})^{\frac{2-n}{2}} e^{i\Phi}$, and we get the equation 
\begin{equation} \label{dbarphieq}
\partial_{\bar{z}} \Phi = -\frac{1}{2}(e_1+i e_r) \cdot \tilde{W}^{\sharp} \quad \text{in } \Psi(\Omega).
\end{equation}
It is now easy to solve \eqref{transporteq_mollified}.

\begin{lemma} \label{lemma:amplitude}
The equation \eqref{transporteq_mollified} has a solution $a \in C^{\infty}(\closure{\Omega})$ which satisfies $\norm{\partial^{\alpha} a}_{L^{\infty}(\Omega)} \leq C_{\alpha} h^{-\sigma\abs{\alpha}}$.
\end{lemma}
\begin{proof}
We first take $\chi \in C^{\infty}_c(\tilde{\Omega})$ with $\chi = 1$ near $\Omega$, and consider 
\begin{equation*}
\partial_{\bar{z}} \Phi = -\frac{1}{2} \tilde{\chi} (e_1+i e_r) \cdot \tilde{W}^{\sharp} \quad \text{for } z \in \mC
\end{equation*}
with $\tilde{\chi} = \chi \circ \Psi^{-1}$. This has the explicit solution 
\begin{equation} \label{phi_explicit}
\Phi(z,\theta) = -\frac{1}{4\pi i} \int_{\mC} \frac{1}{w} \tilde{\chi}(z-w,\theta)(e_1+i e_r) \cdot \tilde{W}^{\sharp}(z-w,\theta) \,d\bar{w} \wedge dw.
\end{equation}
Then $a = ((z-\bar{z})^{\frac{2-n}{2}} e^{i\Phi}) \circ \Psi$ will be a solution of \eqref{transporteq_mollified} in $\Omega$ with the given norm bounds.
\end{proof}

It remains to solve \eqref{correctioneq_mollified}. The $L^2$ norm of the right hand side of \eqref{correctioneq_mollified} is $O(h^{1+\sigma\varepsilon})$ when $\sigma$ is small enough, using the H\"older continuity of $W$. It follows from Proposition \ref{prop:secondinhomog} that when $h$ is small there is an $H^1(\Omega)$ solution $r$ with $\norm{r}_{L^2(\Omega)} = O(h^{\sigma\varepsilon})$, $\norm{\nabla r}_{L^2(\Omega)} = O(h^{-1+\sigma\varepsilon})$.

We collect the results obtained in this argument.
\begin{prop}
  Let $W \in C^{\varepsilon}(\closure{\Omega} ; \mC^n)$, $\nabla \cdot W = 0$,
  and $V \in L^n(\Omega)$. Let $\varphi,\psi$ be defined by
  \eqref{phidef} and \eqref{psidef} respectively. Then for $h$ small
  there is an $H^1(\Omega)$ solution $u =
  e^{-\frac{1}{h}(\varphi+i\psi)}(a+r)$ of the equation $H_{W,V} u = 0$
  in $\Omega$, where $a$ is given in Lemma \ref{lemma:amplitude}.
  Also, one has the norm estimates
\begin{eqnarray*}
 & \norm{\partial^{\alpha} a}_{L^{\infty}(\Omega)} = O(h^{-\sigma\abs{\alpha}}), & \\
 & \norm{r}_{L^2(\Omega)} + \norm{h\nabla r}_{L^2(\Omega)} = O(h^{\sigma\varepsilon}) & 
\end{eqnarray*}
where $\sigma > 0$ is small.
\end{prop}

\section{Determining a magnetic field}

In this section we prove Theorem \ref{thm:magnetic} following the arguments given in \cite{dksu}. Assume that $W_j \in C^{\varepsilon}(\Omega; \mC^n)$, $\nabla \cdot W_j = 0$, $V_j \in L^n(\Omega)$, and that 
\begin{equation} \label{dnmapscoincide}
N_{W_1,V_1} f|_{\tilde{F}} = N_{W_2,V_2} f|_{\tilde{F}} \quad \text{for all } f \in H^{1/2}(\partial \Omega).
\end{equation}
Also assume that $W_j$ are extended as vector fields in $C^{\varepsilon}_c(\tilde{\Omega} ; \mC^n)$.

We start by choosing $u_1$ and $u_2$ to be solutions of $H_{W_1,V_1} u_1 = 0$ and $H_{\bar{W}_2,\bar{V}_2} u_2 = 0$, of the form 
\begin{align*}
u_1 &= e^{\frac{1}{h}(\varphi + i\psi)}(a_1 + r_1), \\
u_2 &= e^{\frac{1}{h}(-\varphi + i\psi)}(a_2 + r_2)
\end{align*}
where $\varphi(x)$ and $\psi(x)$ are defined by \eqref{phidef} and
\eqref{psidef} respectively, $a_j$ satisfies
\begin{eqnarray*}
 & (\nabla(\varphi+i\psi) \cdot D + \nabla(\varphi+i\psi) \cdot W_1^{\sharp} + \frac{1}{2i} \Delta(\varphi+i\psi))a_1 = 0, & \\
 & (\nabla(-\varphi+i\psi) \cdot D + \nabla(-\varphi+i\psi) \cdot \bar{W}_2^{\sharp} + \frac{1}{2i} \Delta(-\varphi+i\psi))a_2 = 0, & 
\end{eqnarray*}
and $\norm{\partial^{\alpha} a_j}_{L^{\infty}(\Omega)} = O(h^{-\sigma\abs{\alpha}})$ for $\sigma > 0$ small, and finally $\norm{r_j} = O(h^{\sigma \varepsilon})$, $\norm{\nabla r_j} = O(h^{-1+\sigma \varepsilon})$ as $h \to 0$.

Take $\tilde{u}_2$ to be the solution to $H_{W_2,V_2} \tilde{u}_2 = 0$ with $\tilde{u}_2 = u_1$ on $\partial \Omega$, and let $u = u_1 - \tilde{u}_2$. We use Lemma \ref{lemma:hdelta_green_identity} for $u$ and $u_2$ to get 
\begin{equation} \label{integralidentity}
(H_{W_2,V_2} u|u_2) = -(\partial_{\nu} u|u_2)_{\partial \Omega}.
\end{equation}
We examine the right hand side of \eqref{integralidentity}.

\begin{lemma} \label{lemma:boundaryterm}
One has $h(\partial_{\nu} u|u_2)_{\partial \Omega} \to 0$ as $h \to 0$.
\end{lemma}
\begin{proof}
Because of \eqref{dnmapscoincide} we have $\partial_{\nu} u = i((W_2-W_1)\cdot \nu)u_1$ on $\tilde{F}$, and so 
\begin{equation} \label{boundarytwoterms}
h(\partial_{\nu} u|u_2)_{\partial \Omega} = ih( ((W_2-W_1)\cdot \nu)u_1 | u_2)_{\tilde{F}} + h(\partial_{\nu} u| u_2)_{\partial \Omega \smallsetminus \tilde{F}}
\end{equation}
The first term satisfies 
\begin{equation*}
h \abs{( ((W_2-W_1)\cdot \nu)u_1 | u_2)_{\tilde{F}}} \lesssim h \norm{a_1+r_1}_{L^2(\partial \Omega)} \norm{a_2+r_2}_{L^2(\partial \Omega)}.
\end{equation*}
We have $\norm{a_j}_{L^2(\partial \Omega)} \lesssim \norm{a_j}_{L^{\infty}(\partial \Omega)} = O(1)$ and $\norm{r_j}_{L^2(\partial \Omega)} \lesssim \norm{r_j}_{H^{1/2+\delta}(\Omega)} = O(h^{-1/2+\sigma\varepsilon-\delta})$ for any $\delta > 0$, which gives 
\begin{equation} \label{ajrjestimate}
\norm{a_j + r_j}_{L^2(\partial \Omega)} = O(h^{-1/2+\delta})
\end{equation}
for some new $\delta > 0$. This shows that the first term in \eqref{boundarytwoterms} vanishes as $h \to 0$.

For the other term we compute 
\begin{equation*}
\abs{h(\partial_{\nu} u|u_2)_{\partial \Omega \smallsetminus \tilde{F}}} \leq h \norm{e^{-\frac{\varphi}{h}} \partial_{\nu} u}_{L^2(\partial \Omega \smallsetminus \tilde{F})} \norm{a_2 + r_2}_{L^2(\partial \Omega)}.
\end{equation*}
To estimate the normal derivative we use the Carleman estimate of Proposition \ref{prop:carleman}, for the weight $-\varphi$, in the form
\begin{multline*}
\sqrt{h} \norm{\sqrt{\partial_{\nu} \varphi} e^{-\frac{\varphi}{h}} \partial_{\nu} u}_{L^2(\partial \Omega_{+})} + \norm{e^{-\frac{\varphi}{h}} u} + \norm{e^{-\frac{\varphi}{h}} h \nabla u} \\
 \lesssim h \norm{e^{-\frac{\varphi}{h}} H_{W_2,V_2} u} + \sqrt{h} \norm{\sqrt{-\partial_{\nu} \varphi} e^{-\frac{\varphi}{h}} \partial_{\nu} u}_{L^2(\partial \Omega_{-})}.
\end{multline*}
The estimate applies for this $u$ since $u \in H_{\Delta}(\Omega)$ and $u|_{\partial \Omega} = 0$. Since $e^{-\frac{\varphi}{h}} \partial_{\nu} u = ie^{i\frac{\psi}{h}}((W_2-W_1)\cdot \nu)(a_1+r_1)$ on $\tilde{F}$, and since $\partial \Omega_{-} \subseteq \tilde{F}$, it follows from \eqref{ajrjestimate} that 
\begin{align*}
\norm{e^{-\frac{\varphi}{h}} \partial_{\nu} u}_{L^2(\partial \Omega \smallsetminus \tilde{F})} &\lesssim \norm{\sqrt{\partial_{\nu} \varphi} e^{-\frac{\varphi}{h}} \partial_{\nu} u}_{L^2(\partial \Omega_{+})} \\
 &\lesssim \sqrt{h} \norm{e^{-\frac{\varphi}{h}} H_{W_2,V_2} u} + O(h^{-1/2+\delta}).
\end{align*}
Using that $u_1$ and $\tilde{u}_2$ are solutions and $\nabla \cdot W_j = 0$, we obtain 
\begin{equation*}
H_{W_2,V_2} u = H_{W_2,V_2} u_1 = 2(W_2-W_1) \cdot Du_1 + (W_2^2 - W_1^2 + V_2 - V_1) u_1.
\end{equation*}
The explicit form for $u_1$ implies 
\begin{multline} \label{hw2v2u}
H_{W_2,V_2} u = e^{\frac{1}{h}(\varphi+i\psi)}(2ih^{-1}[(W_1-W_2) \cdot \nabla (\varphi+i\psi)] (a_1 + r_1) \\
 + 2(W_2-W_1) \cdot D(a_1+r_1) + (W_2^2 - W_1^2 + V_2 - V_1)(a_1 + r_1))
\end{multline}
which shows that $\norm{e^{-\frac{\varphi}{h}} H_{W_2,V_2} u} = O(h^{-1})$, since the $L^{\frac{2n}{n-2}}$ norm of $a_1+r_1$ can be estimated by the $H^1$ norm. Collecting these estimates gives 
\begin{equation*}
h(\partial_{\nu} u|u_2)_{\partial \Omega \smallsetminus \tilde{F}} = O(h^{\delta})
\end{equation*}
as $h \to 0$, which concludes the proof.
\end{proof}

From \eqref{hw2v2u} we obtain 
\begin{equation*}
h(H_{W_2,V_2}u|u_2) = \int_{\Omega} 2i(W_1-W_2) \cdot \nabla(\varphi+i\psi) a_1 \bar{a}_2 \,dx + o(1).
\end{equation*}
Then \eqref{integralidentity} and Lemma \ref{lemma:boundaryterm} imply 
\begin{equation} \label{firsttransformomega}
\int_{\Omega} \nabla(\varphi+i\psi) \cdot (W_1-W_2) a \,dx = 0
\end{equation}
where $a = \lim_{h \to 0} a_1 \bar{a}_2$. We will now switch to the complex notation as in Section \ref{sec:construction}.

Choose coordinates so that $x_0 = 0$ and $\omega = e_1$, and write $x = (x_1,x')$ where $x' = r\theta$, $r > 0$, and $\theta \in S^{n-2}$. Let $P_{\theta}$ be the two-dimensional plane consisting of points $(x_1,r\theta)$ for $\theta$ fixed, and write $\Omega_{\theta} = \Omega \cap P_{\theta}$. We also use the complex variable $z = x_1 + ir$, which identifies $P_{\theta}$ with $\mC$.

In the coordinates $(z,\theta) = \Psi(x)$, Lemma \ref{lemma:amplitude} shows that $a$ in \eqref{firsttransformomega} is given by $a \circ \Psi^{-1} = (z-\bar{z})^{2-n} e^{i\Phi}$, where $\Phi$ satisfies 
\begin{equation} \label{phitransport}
\partial_{\bar{z}} \Phi = -\frac{1}{2}(e_1+i e_r) \cdot (\tilde{W}_1 - \tilde{W}_2) \quad \text{in } \Psi(\Omega).
\end{equation}
We have written $\tilde{W}_j = W_j \circ \Psi^{-1}$. Note that since $\Phi$ is of the form \eqref{phi_explicit}, $W_j \in C^{\varepsilon}_c$ implies that $\Phi \in C^{\varepsilon}$ and further $\Phi(\,\cdot\,,\theta) \in C^{1+\varepsilon}(\mC)$. The latter follows since $\partial_{\bar{z}} \Phi \in C^{\varepsilon}$ and $\partial_{z} \Phi \in \partial_{z} \partial_{\bar{z}}^{-1} C^{\varepsilon}_c$, and the singular integral operator $\partial_{z} \partial_{\bar{z}}^{-1}$ is bounded on H\"older spaces.

Now \eqref{firsttransformomega} becomes 
\begin{equation} \label{firsttransform}
\int_{S^{n-2}} \Big( \int_{\Omega_{\theta}} \frac{1}{z} (e_1+ie_r) \cdot (\tilde{W}_1-\tilde{W}_2) e^{i\Phi} \,d\bar{z} \wedge dz \Big) \,d\theta = 0.
\end{equation}
We need a slightly more general result. From the transport equation for $a_1$ we see that in the definition of $u_1$ we may replace $a_1$ by $a_1 g_1$, where $\nabla(\varphi+i\psi) \cdot \nabla g_1 = 0$ and $g_1 \in W^{2,\infty}(\Omega)$ is independent of $h$. It follows that \eqref{firsttransform} holds with $e^{i\Phi}$ replaced by $e^{i\Phi} g_1$. Choosing $g_1 = z g(z) \tilde{g}(\theta)$ with $\partial_{\bar{z}} g = 0$ and $\tilde{g}$ smooth, and by varying $\tilde{g}$, we see that for almost every $\theta \in S^{n-2}$ one has 
\begin{equation} \label{secondtransform}
\int_{\Omega_{\theta}} (e_1+ie_r) \cdot (\tilde{W}_1-\tilde{W}_2) e^{i\Phi} g(z) \,d\bar{z} \wedge dz = 0.
\end{equation}
The argument that \eqref{secondtransform} implies $dW_1 = dW_2$ now proceeds as in \cite{dksu}. However, this was written in \cite{dksu} quite briefly, especially in the case where $\Omega_{\theta}$ has nontrivial topology. Therefore, we will give a rather detailed argument for determining the coefficients.

It will be enough to restrict to $3$-dimensional subspaces $L$ of $\mR^n$ which consist of the points $(x_1,r\theta)$ where $\theta \in S^{n-2}$ varies on a fixed two-plane. The fact that \eqref{secondtransform} holds in each $L$ will imply $dW_1 = dW_2$ in $\Omega$, see \cite{dksu}. Thus we may assume $n = 3$ in what follows, and $\theta \in S^1$ may be identified with an angle, using $\angle e^{i\alpha} = \alpha$ when $0 < \alpha < \pi$.

We begin by looking at what kind of topology $\Omega_{\theta}$ can have when $\theta$ varies. For this we use the map $\Theta: \partial \Omega \to \mR, x \mapsto \angle \frac{x'}{\abs{x'}}$. The main point is that possibly after a small change of coordinates, $\Theta$ is a Morse function.

\begin{lemma} \label{lemma:morse}
For almost every choice $(\omega,x_0) \in T S^2$, $\Theta$ is a Morse function and there are finitely many critical values $\theta_0 < \theta_1 < \ldots < \theta_N$ of $\Theta$. Also, $\Theta(\partial \Omega) = [\theta_0,\theta_N]$, and if $\theta \in (\theta_j,\theta_{j+1})$ then 
\begin{enumerate}
\item[(a)] 
$\partial \Omega$ and $P_{\theta}$ intersect transversally, or equivalently for any $x \in \partial \Omega \cap P_{\theta}$ the vector $\nu(x)$ is not orthogonal to $P_{\theta}$,
\item[(b)]
$\Omega_{\theta}$ is a bounded open set with smooth boundary in $P_{\theta}$, and $\partial \Omega_{\theta} = \partial \Omega \cap P_{\theta}$,
\item[(c)]
there is a diffeomorphism $T:(\theta_j,\theta_{j+1}) \times \partial \Omega_{\theta} \to \Theta^{-1}((\theta_j,\theta_{j+1}))$, and if $T_t(x) = T(t,x)$ then $T_t$ is a diffeomorphism of $\partial \Omega_{\theta}$ and $\partial \Omega_t$ with $T_{\theta} = \text{Id}_{\partial \Omega_{\theta}}$.
\end{enumerate}
\end{lemma}
\begin{proof}
We see that $x \in \partial \Omega$ is a degenerate critical point of $\Theta$ if and only if $N(x) \cdot \omega = N(x) \cdot (x-x_0) = 0$ and $K(x) = 0$, where $N$ is the Gauss map and $K$ is the Gaussian curvature of $\partial \Omega$. Thus, $\Theta$ is Morse provided that $(\omega,x_0)$ is not a critical value of the map $F: \Omega(\partial \Omega) \times \mR \to T S^2$, 
\begin{equation*}
F((x,\omega),t) = (\omega, (x \cdot N(x)) N(x) + t(\omega \times N(x))),
\end{equation*}
where $\Omega (\partial \Omega) \subseteq T(\partial \Omega)$ is the unit sphere bundle. Sard's lemma shows that $\Theta$ is Morse for almost every $(\omega,x_0)$. Then $\Theta$ has finitely many critical values, and if $\theta$ is not a critical value then $\partial \Omega$ and $P_{\theta}$ intersect transversally.

The other condition in (a) is just another way of stating the transversality, and this latter condition implies that $\partial \Omega_{\theta} = \partial \Omega \cap P_{\theta}$. Then, $\partial \Omega_{\theta}$ is a smooth $1$-manifold and $\Omega_{\theta}$ has smooth boundary in $P_{\theta}$. Part (c) follows from the fact that there is a gradient-like vector field for $\Theta$ \cite{milnor}, and we obtain the diffeomorphism $T$ by considering flows of this vector field which originate from $\partial \Omega_{\theta}$.
\end{proof}

Note that $\partial \Omega_{\theta}$ can have many components, but if $\theta$ stays away from the critical values then the number of components stays fixed and the components vary smoothly with $\theta$. Unless stated otherwise, we will assume that the coordinates are chosen as in Lemma \ref{lemma:morse} and $\theta$ is not a critical value.

\begin{lemma} \label{lemma:curlidentical}
$dW_1 = dW_2$ in $\Omega$.
\end{lemma}
\begin{proof}
Our starting point is \eqref{secondtransform}. Because of the factor $e^{i\Phi}$ this may be considered as a nonlinear Radon transform, evaluated at the plane $P_{\theta}$. Most of the work below is to show that \eqref{secondtransform} remains true with $e^{i\Phi}$ replaced by $1$, which corresponds to the usual Radon transform.

The proof is in several steps. For the complex analysis terminology see \cite{rudin}.

\medskip
\noindent \emph{Step 1}.
Using the equation \eqref{phitransport} and integrating by parts, we obtain from \eqref{secondtransform} the orthogonality condition 
\begin{equation} \label{thirdtransform}
\int_{\partial \Omega_{\theta}} e^{i\Phi} g(z) \,dz = 0
\end{equation}
for any holomorphic $g \in W^{2,\infty}(\Omega_{\theta})$. Since $P_{\theta} \smallsetminus \Omega_{\theta}$ has finitely many components, an approximation argument (using a version of Mergelyan's theorem) implies that \eqref{thirdtransform} holds in fact for any holomorphic $g \in C(\closure{\Omega}_{\theta})$.

\medskip
\noindent \emph{Step 2}.
The condition \eqref{thirdtransform} is equivalent with saying that $e^{i\Phi}|_{\partial \Omega_{\theta}}$ is the boundary value of a holomorphic function $F \in C(\closure{\Omega}_{\theta})$. Indeed, write $f = e^{i\Phi}$, and let $F$ be the Cauchy integral of $f|_{\partial \Omega_{\theta}}$. From \eqref{thirdtransform} we see that $F = 0$ outside $\Omega_{\theta}$, and the Plemelj-Sokhotski formula then implies that $F \in C(\closure{\Omega}_{\theta})$ with $F|_{\partial \Omega_{\theta}} = f$.

\medskip
\noindent \emph{Step 3}.
We would like to show that $F$ is nonvanishing and has a holomorphic logarithm in $\Omega_{\theta}$. We first claim that if $\gamma$ is a closed curve in $\partial \Omega_{\theta}$, then 
\begin{equation} \label{firstindexzero}
\int_{F \circ \gamma} \frac{1}{z} \,dz = \int_{f \circ \gamma} \frac{1}{z} \,dz = 0.
\end{equation}
To show this we write $f_s = e^{is\Phi}$ for $0 \leq s \leq 1$. We obtain that $f \circ \gamma$ is homotopic to $f_0 \circ \gamma = \{1\}$ in $\mC \smallsetminus \{0\}$, and the claim follows.

\medskip
\noindent \emph{Step 4}.
We can use \eqref{firstindexzero} and the argument principle to conclude that $F$ is nonvanishing, also when $\Omega_{\theta}$ is not simply connected. Let $\gamma$ be an oriented parametrization of $\partial \Omega_{\theta}$ as a sum of simple closed curves. It follows that $\text{Ind}_{\gamma}(\alpha) = \int_{\partial \Omega_{\theta}} \frac{1}{z-\alpha} \,dz = 0$ for all $\alpha \notin \Omega_{\theta}$. Also, $F \neq 0$ on $\partial \Omega_{\theta}$, so there are only finitely many zeros in $\Omega_{\theta}$. One may now use the residue theorem and argue in the usual way that $\frac{1}{2\pi i} \int_{F \circ \gamma} \frac{1}{z} \,dz = 0$ is the number of zeros of $F$ in $\Omega_{\theta}$.

\medskip
\noindent \emph{Step 5}.
Next we will show that $F$ has a holomorphic logarithm in $\Omega_{\theta}$. Again, this would be immediate in a simply connected domain. In the general case, $F$ has a holomorphic logarithm provided that 
\begin{equation} \label{logarithmcondition}
\int_{\gamma} \frac{F'}{F} \,dz = 0
\end{equation}
for any closed curve $\gamma$ in $\Omega_{\theta}$. To show \eqref{logarithmcondition} for given $\gamma$, we take $C_j$ to be the finitely many components of $P_{\theta} \smallsetminus \Omega_{\theta}$, and for each $j$ we let $\gamma_j$ be a cycle corresponding to the oriented boundary of $C_j$. Then $\gamma_j$ is contained in $\partial \Omega_{\theta}$, the index of $\gamma_j$ on $C_j$ is $1$, and the index of $\gamma_j$ on any other component $C_k$ is $0$. Thus, using the $\gamma_j$, we can construct a cycle $\tilde{\gamma}$ contained in $\partial \Omega_{\theta}$ so that the index of $\tilde{\gamma}$ is equal to the index of $\gamma$ at each point outside of $\Omega_{\theta}$. It follows that 
\begin{equation*}
\int_{\gamma} \frac{F'}{F} \,dz = \int_{\tilde{\gamma}} \frac{F'}{F} \,dz = \int_{F \circ \tilde{\gamma}} \frac{1}{z} \,dz = 0
\end{equation*}
using \eqref{firstindexzero}. This shows \eqref{logarithmcondition}, and we obtain a holomorphic logarithm $G \in C(\closure{\Omega}_{\theta})$ with $F = e^G$ by fixing a point $z_0$ in each component of $\Omega_{\theta}$, and by taking in this component 
\begin{equation} \label{gdef}
G(z) = \int_{z_0}^z \frac{F'}{F} \,dw + c_0
\end{equation}
where the integral is over any curve connecting $z_0$ to $z$, and $e^{c_0} = F(z_0)$.

\medskip
\noindent \emph{Step 6}.
Since $e^{G-i\Phi} = 1$ on $\partial \Omega_{\theta}$, we get $G|_{\partial \Omega_{\theta}} = i\Phi + v$ where $v$ is constant on each component of $\partial \Omega_{\theta}$. In fact, $v$ is equal to $2\pi i m$ for some $m \in \mZ$ on each of these components. This shows that 
\begin{equation} \label{gboundaryintegral}
\int_{\partial \Omega_{\theta}} G \,dz = \int_{\partial \Omega_{\theta}} i\Phi \,dz.
\end{equation}

\medskip 
\noindent \emph{Step 7}.
We now return to \eqref{thirdtransform} and take $g = \frac{G}{e^G} \in C(\closure{\Omega}_{\theta})$. Then \eqref{gboundaryintegral} gives 
\begin{equation*}
\int_{\partial \Omega_{\theta}} \Phi \,dz = 0.
\end{equation*}
Integrating by parts and using \eqref{phitransport} yields 
\begin{equation} \label{lasttransform}
\int_{\Omega_{\theta}} (e_1+ie_r) \cdot (\tilde{W}_1-\tilde{W}_2) \,d\bar{z} \wedge dz = 0.
\end{equation}

\medskip
\noindent \emph{Step 8}.
We need to show that \eqref{lasttransform} is valid with $e_r$ replaced by $-e_r$. If $W_j$ were real valued this would follow just by taking complex conjugates. If the $W_j$ are complex, we can go back to the beginning of this section and repeat the construction of solutions, with $\psi$ replaced by $-\psi$. In this way, instead of \eqref{firsttransformomega} we arrive at 
\begin{equation*}
\int_{\Omega} \nabla (\varphi+i\psi) \cdot (\bar{W}_1 - \bar{W}_2) a \,dx = 0
\end{equation*}
where $a \circ \Psi^{-1} = (z-\bar{z})^{2-n} e^{i\Phi}$, and $\Phi$ satisfies 
\begin{equation*}
\partial_{\bar{z}} \Phi = \frac{1}{2}(e_1+i e_r) \cdot (\widetilde{\bar{W}}_1 - \widetilde{\bar{W}}_2) \quad \text{in } \Psi(\Omega).
\end{equation*}
Going through steps 1 to 8 above, we obtain \eqref{lasttransform} with $W_j$ replaced by $\bar{W}_j$, and taking conjugates gives 
\begin{equation} \label{lasttransformconjugate}
\int_{\Omega_{\theta}} (e_1-ie_r) \cdot (\tilde{W}_1-\tilde{W}_2) \,d\bar{z} \wedge dz = 0.
\end{equation}
 
\medskip
\noindent \emph{Step 9}.
If $\xi \in P_{\theta}$ then $\xi$ is a linear combination of $e_1$ and $e_r$, and we get from \eqref{lasttransform} and \eqref{lasttransformconjugate} that 
\begin{equation*}
\int_{\Omega_{\theta}} \xi \cdot (\tilde{W}_1-\tilde{W}_2) \,d\bar{z} \wedge dz = 0.
\end{equation*}
Returning to the $x$ coordinates, this gives 
\begin{equation} \label{radontransform}
\int_{(x_0+P) \cap \Omega} \xi \cdot (W_1 - W_2) \,dS = 0
\end{equation}
for all two-planes $P$ passing through $e_1$, and all $\xi \in P$.

\medskip
\noindent \emph{Step 10}.
The left hand side of \eqref{radontransform} is related to the Radon transform of an expression involving $\curl\,(W_1-W_2)$. Varying $x_0$ in a small neighborhood of $0$ and varying $\omega$ in a conic neighborhood of $e_1$, one obtains $\curl\,(W_1-W_2) = 0$ from the arguments in \cite[Lemma 5.2]{dksu}.
\end{proof}

After showing that $dW_1 = dW_2$, the final step is to show that $V_1 = V_2$. This also follows along the lines of \cite{dksu}, but we need to give more details to account for the nontrivial topology of the $\Omega_{\theta}$.

\begin{lemma}
$V_1 = V_2$ in $\Omega$.
\end{lemma}
\begin{proof}
The proof is again in several steps.

\medskip
\noindent \emph{Step 1}.
Since $\curl(W_1-W_2) = 0$ and $\Omega$ is simply connected, we have $W_2 - W_1 = \nabla p$ where $p \in C^{1+\varepsilon}(\closure{\Omega})$.

\medskip
\noindent \emph{Step 2}.
The assumption that $N_{W_1,V_1} f|_{\tilde{F}} = N_{W_2,V_2} f|_{\tilde{F}}$ for all $f \in H^{1/2}(\partial \Omega)$, together with boundary determination results for the magnetic Schr\"odinger equation \cite{brownsalo}, implies that $(W_1)_{\tan} = (W_2)_{\tan}$ on $\tilde{F}$. Note that the results in \cite{brownsalo} remain valid for $V \in L^n(\Omega)$. It follows that the tangential derivatives of $p$ vanish on $\tilde{F}$, which shows that $p$ is constant on each component of $\tilde{F}$.

\medskip
\noindent \emph{Step 3}.
We use again the complex notation, with $x_0 = 0$ and $\omega = e_1$. Inserting $W_2 - W_1 = \nabla p$ in \eqref{phitransport}, we see that $\Phi = \tilde{p}$ will be a solution of \eqref{phitransport} if $\tilde{p} = p \circ \Psi^{-1}$. Repeating the arguments in the proof of Lemma \ref{lemma:curlidentical} shows that $e^{i\tilde{p}}|_{\partial \Omega_{\theta}}$ is the boundary value of a nonvanishing holomorphic function, and taking a holomorphic logarithm gives $i\tilde{p} = G_{\theta} + v_{\theta}$ on $\partial \Omega_{\theta}$, where $G_{\theta} \in C(\closure{\Omega}_{\theta})$ is holomorphic and $v_{\theta}$ is equal to some $2\pi i m$, $m \in \mZ$, on each component of $\partial \Omega_{\theta}$.

\medskip
\noindent \emph{Step 4}.
We claim that for each component $C_{\theta}$ of $\Omega_{\theta}$, there is a point $y_0 \in \partial C_{\theta}$ with $y_0 \in \tilde{F}$. In fact, we may take $y_0$ to be a point which minimizes $\abs{y}$ among $y \in \closure{C}_{\theta}$. Then $y_0 \in \partial C_{\theta}$ and $(1-t)y_0 \notin \closure{C}_{\theta}$ for $t > 0$. It is enough to show that $y_0 \cdot \nu(y_0) \leq 0$, since then $y_0 \in \tilde{F}$. But if one had $y_0 \cdot \nu(y_0) > 0$, then the fact that for any $c > 0$ there is a truncated cone 
\begin{equation*}
\{ y \in \mR^n \colon (y-y_0) \cdot \nu(y_0) < -c\abs{y-y_0}, \abs{y-y_0} < \delta \} \cap \Omega_{\theta}
\end{equation*}
contained in $C_{\theta}$, would imply that $(1-t)y_0 \in C_{\theta}$ for some $t > 0$. This is a contradiction.

\medskip
\noindent \emph{Step 5}.
From Steps 3 and 4 we see that $G_{\theta} \in C(\closure{\Omega}_{\theta})$ is a holomorphic function in $\Omega_{\theta}$, and for any component $C_{\theta}$ of $\Omega_{\theta}$ there is an open set of $\partial C_{\theta}$ in which $G_{\theta}$ is constant. This implies that $G_{\theta}$ is constant on each component of $\closure{\Omega}_{\theta}$, hence also on each component of $\partial \Omega_{\theta}$.

\medskip
\noindent \emph{Step 6}.
We also need that $G_{\theta}|_{\partial \Omega_{\theta}}$ can be made to vary continuously with $\theta$, as long as $\theta$ stays away from the critical values of $\Theta$. To do this we use Lemma \ref{lemma:morse} (c) which implies that $\partial \Omega_{\theta}$ varies continuously with $\theta$, and so does $F$ in \eqref{gdef}. Using the regularity properties of $\Phi = \tilde{p}$, and choosing the point $z_0$ in \eqref{gdef} to lie in $\partial \Omega_{\theta}$ with continuous dependence on $\theta$, it can be checked that $G_{\theta}|_{\partial \Omega_{\theta}}$ varies continuously with $\theta$.

\medskip
\noindent \emph{Step 7}.
We now restrict to the set $\Theta^{-1}((\theta_j,\theta_{j+1}))$, and claim that $p$ is locally constant in this set (which implies that $p$ is constant on the components of this set, since $p$ is continuous). To prove this, we first recall that $i\tilde{p} = G_{\theta} + v_{\theta}$ on $\partial \Omega_{\theta}$, where $\tilde{p}$ and $G_{\theta}$ vary continuously with $\theta$. Then also $v_{\theta}$ varies continuously with $\theta$, and since $\frac{1}{2\pi i} v_{\theta}$ is integer valued we see that $v_{\theta}$ is locally constant.

\medskip
\noindent If $x \in \Theta^{-1}((\theta_j,\theta_{j+1}))$ and $\Theta(x) = \theta$, then $x$ lies in some component $C_{\theta}$ of $\closure{\Omega}_{\theta}$. By Step 4 there is some point $y$ in $C_{\theta} \cap \tilde{F}$, and Step 2 shows that $p$ is constant in $B(y,r) \cap \partial \Omega$ for some $r > 0$. The same then applies to $G_{\theta}$ near $y$. Since the components of $\partial \Omega_{\theta}$ vary continuously, and since $G_{\theta}$ is constant on each such component, we see that $G_{\theta}$ is constant near $x$ when $\theta$ varies. This shows the claim.

\medskip
\noindent \emph{Step 8}.
To show that $p$ is locally constant on $\partial \Omega$, it remains to check that this is true near each $x \in \partial \Omega \cap P_{\theta}$ where $\theta$ is a critical value. If $\nu(x)$ is orthogonal to $P_{\theta}$, then $x \in \tilde{F}$ and this follows by Step 2. If $\nu(x)$ is not orthogonal to $P_{\theta}$, then for some $r > 0$ and $j$, both the sets $B(x,r) \cap \Theta^{-1}((\theta_j,\theta_{j+1}))$ and $B(x,r) \cap \Theta^{-1}((\theta_{j-1},\theta_j))$ are nonempty and connected. Since $p$ is constant on both these sets and continuous, $p$ must be constant also near $x$.

\medskip
\noindent \emph{Step 9}.
The preceding step implies that $p$ is constant on the components of $\partial \Omega$. This shows that $(W_1)_{\text{tan}} = (W_2)_{\text{tan}}$ on $\partial \Omega$. Recalling that $\partial \Omega$ is connected, we may substract a constant from $p$ to obtain $p|_{\partial \Omega} = 0$. The assumption, gauge invariance, and Step 1 then imply 
\begin{equation*}
N_{W_1,V_1}|_{\tilde{F}} = N_{W_2,V_2}|_{\tilde{F}} = N_{W_2-\nabla p,V_2}|_{\tilde{F}} = N_{W_1,V_2}|_{\tilde{F}}.
\end{equation*}
Consequently, we may assume $W_1 = W_2$ in the arguments in this section. Going through the proof of Lemma \ref{lemma:boundaryterm}, and using the assumption $V_j \in L^{\infty}(\Omega)$, we may take limits as $h \to 0$ in \eqref{integralidentity} to obtain 
\begin{equation} \label{firstvdiff}
\int_{\Omega} (V_1-V_2) a \,dx = 0
\end{equation}
where $a \circ \Psi^{-1} = (z-\bar{z})^{2-n}$.

\medskip
\noindent \emph{Step 10}.
We argue as in Lemma \ref{lemma:curlidentical} and replace $a$ in \eqref{firstvdiff} by $ag$ where $\nabla(\varphi+i\psi) \cdot \nabla g = 0$. Moving to the variables $(z,\theta)$ and taking $g$ to be a function of $\theta$, and by varying $g$, we obtain 
\begin{equation} \label{secondvdiff}
\int_{\Omega_{\theta}} (V_1-V_2)(\Psi^{-1}(z,\theta)) \,d\bar{z} \wedge dz = 0 
\end{equation}
for almost every $\theta$. Here we use the fact that $V_1-V_2$ is in $L^1(\mR^n)$, so the restriction $V_1-V_2|_{P_{\theta}}$ is integrable on $P_{\theta}$ for almost every $\theta$.

\medskip
\noindent In the $x$ coordinates, \eqref{secondvdiff} reads 
\begin{equation}
\int_{(x_0+P) \cap \Omega} (V_1 - V_2) \,dS = 0.
\end{equation}
This is valid for almost every two-plane $P$ passing through $\omega = e_1$. Varying $x_0$ and $\omega$ slightly, the Radon transform arguments in \cite{dksu} give $V_1 = V_2$.
\end{proof}

\section{Determining a convection term}

Finally, we prove Theorem \ref{thm:convection} by reducing the inverse problem for the convection equation to the corresponding problem for the magnetic Schr\"odinger equation.

\begin{proof}
(of Theorem \ref{thm:convection}) We observe that $H_{iW_j,q_j} = -\Delta + 2W_j \cdot \nabla$ where $q_j = W_j^2 - \nabla \cdot W_j$. Consequently, for $f \in H^{1/2}(\partial \Omega)$ we have the relation $N_{iW_j,q_j} f = N_{W_j} f - (W_j \cdot \nu) f$ between DN maps, and 
\begin{equation*}
(N_{iW_1,q_1} - N_{iW_2,q_2})f|_{\tilde{F}} = -((W_1-W_2) \cdot \nu)f|_{\tilde{F}}
\end{equation*}
using the assumption that $N_{W_1} f = N_{W_2} f$ on $\tilde{F}$. By \cite{brownsalo}, this assumption also implies that $W_1 = W_2$ on $\tilde{F}$, and we obtain $(N_{iW_1,q_1} - N_{iW_2,q_2})f|_{\tilde{F}} = 0$.

Now Theorem \ref{thm:magnetic} shows that $dW_1 = dW_2$ and $q_1 = q_2$ in $\Omega$. Since $\Omega$ is simply connected there is $p \in C^{1+\varepsilon}(\closure{\Omega})$ with $W_2 = W_1 + \nabla p$. Theorem \ref{thm:magnetic} also shows that $(W_1)_{\mathrm{tan}} = (W_2)_{\mathrm{tan}}$ on $\partial \Omega$, so $p$ is constant on $\partial \Omega$ which was assumed to be connected. By substracting the constant we may assume that $p|_{\partial \Omega} = 0$. The condition $q_1 = q_2$ then implies that 
\begin{equation*}
-\Delta p + (2 W_1 + \nabla p) \cdot \nabla p = 0 \qquad \text{in } \Omega.
\end{equation*}
Since everything is real valued, the maximum principle shows that the only solution with zero boundary values is $p \equiv 0$. It follows that $W_1 = W_2$.
\end{proof}

\section{Acknowledgements}

K.K.~was supported by the Carlsberg Foundation, and M.S.~was supported by the Academy of Finland. Part of this research was carried out while the authors were visiting the University of Washington, and while the second author was visiting Aalborg University. Both authors would like to thank Gunther Uhlmann for his generous support.

\addcontentsline{toc}{chapter}{Bibliography}
\bibliography{magnetic_local_nonsmooth}
\bibliographystyle{hamsplain}

\end{document}